%% file: main.tex
\documentclass[11pt,a4paper,reqno]{amsart}

\usepackage{hyperref}
\hypersetup{ 
  pdfauthor   = {Jean-Marc Couveignes, Tony Ezome and Reynald Lercier}
  pdftitle    = {A Faster Pseudo-Primality Test}
  pdfsubject  = {11Y11},
  pdfkeywords = {Primality, Ring theory,  Galois theory, Probabilistic algorithms},
  backref=true, pagebackref=true, hyperindex=true, colorlinks=true,
  breaklinks=true, urlcolor=blue, linkcolor=blue, citecolor=blue, bookmarks=true,
  bookmarksopen=true}

\newcommand{\secstar}[1]{\subsection*{#1}}
\newcommand{\subsecstar}[1]{\subsubsection*{#1}}

\newenvironment{myproof}[1][\myproofname]{\par
  \normalfont \topsep6pt\relax
  \trivlist
\item[\hskip\labelsep
  \itshape
  #1.]\ignorespaces
}{%
  \endtrivlist\hfill$\square$
}
\providecommand{\myproofname}{Proof}

\usepackage{multirow}

\usepackage{xy,amsmath,amssymb,url,a4wide,amsbsy} 
\usepackage{graphicx,color}

\xyoption{all}

\usepackage[english]{babel}

\let\set\mathbb
\def\<<{\leavevmode
  \raise0.28ex\hbox{$\scriptscriptstyle\langle\!\langle$}\nobreak
  \hskip -.6pt plus.3pt minus.2pt\,}
\def\>>{\,\nobreak\hskip -.6pt plus.3pt minus.2pt
  \raise0.28ex\hbox{$\scriptscriptstyle\rangle\!\rangle$}}

\def\dc{{d_{\rm cyc}}}
\def\xc{u}
\def\sigmac{{\sigma_{\rm cyc}}}
\def\Rc{{R_{\rm cyc}}}
\def\Rcs{{R^{*}_{\rm cyc}}}
\def\dk{{d_{\rm kum}}}
\def\max{{\mathop{{\rm max}}\nolimits }}

\def\rmP{{\mathrm P}}
\def\mprime{{\mathrm{prime}}}
\def\composite{{\mathrm{composite}}}
\def\rmMR{{\mathrm{MR}}}

\def\Log{\mathop{\rm{Log}}\nolimits }
\def\Exp{\mathop{\rm{Exp}}\nolimits }

\def\bK{{\bf K  }}
\def\bM{{\bf M  }}

\def\bL{{\bf L  }}

\def\bL{{\bf L  }}

\def\RR{{\set R}}

\def\Id{{\rm{Id}}}
\def\UU{{\set U}}
\def\VV{{\set V}}
\def\ZZ{{\set Z}}

\def\Ngot{{\mathfrak N}}

\def\cA{{\mathcal A}}
\def\cB{{\mathcal B}}

\def\cG{{\mathcal G}}

\def\pgot{{\mathfrak p}}

\newtheorem{theorem}{Theorem}

\usepackage{enumitem}

\title{A Faster Pseudo-Primality Test} 

\author{Jean-Marc Couveignes}%
\address{Institut de Math\'ematiques de Bordeaux, Universit\'e  Bordeaux I et CNRS, 351 cours de la Lib{\'e}ration, 33405 Talence cedex,  France.}
\address{INRIA Bordeaux Sud-Ouest, projet LFANT.}
\email{jean-marc.couveignes@math.u-bordeaux1.fr}
\author{Tony Ezome}
\address{Universit{\'e} des Sciences et Techniques de Masuku,
Facult{\'e} des Sciences, D{\'e}partement de math{\'e}matiques et informatique,
BP 943 Franceville, Gabon.}%
\email{latonyo2000@yahoo.fr}

\author{Reynald Lercier}%
\address{%
  \textsc{DGA MI}, %
  La Roche Marguerite, %
  35174 Bruz, %
  France. %
}
\address{%
  Institut de recherche math\'ematique de Rennes, %
  Universit\'e de Rennes 1, %
  Campus de Beaulieu, %
  35042 Rennes, %
  France. %
} 
\email{reynald.lercier@m4x.org}

\thanks{Research supported by the ``Direction  G{\'e}n{\'e}rale de 
    l'Armement'', by the ``Agence Nationale de la Recherche'' (projects
    ALGOL and CHIC), and by ``INRIA'' (project LFANT).}

\date{\today}

\begin{document}

\begin{abstract}
We propose a  pseudo-primality test  using cyclic
extensions
of $\ZZ/n\ZZ$. 
For every positive integer  $k \leqslant \log n$, this   test
achieves the    security of $k$
Miller-Rabin tests at   the  cost  of $k^{1/2+o(1)}$
Miller-Rabin tests. 
\end{abstract}


\maketitle

\section{Introduction}\label{sec:intro}

\secstar{Pseudo-primality tests} 
The most commonly used algorithm for prime detection is the so called
Miller-Rabin test.  It is a Monte Carlo probabilistic test of compositeness,
also called a {\it pseudo-primality test} (see Papadimitrou's book \cite[page
254]{PAPA} for the definition of a Monte Carlo algorithm).  A pseudo-primality
test is a process based on a mathematical statement, the {\it compositeness
  criterion}, which gives a forecast (prime or composite) about a given
integer $n$.  From the compositeness criterion, one constructs for every odd
integer $n$, a finite set $W_n$ of {\it witnesses}, and a map
\[\mathrm{P}_n : W_n \rightarrow \{\composite, \mprime \}\]
which provides  information about the compositeness of $n$ from 
witnesses $x$ in $W_n$.
When $n$ is prime $\rmP_n(x)=\mprime$  for every witness
$x$ in $W_n$. So there are  only {\it good witnesses} in that case. 
If $n$ is composite, $x$ is a witness in $W_n$, and $\rmP_n(x)=\mprime$
we say that $x$ is a {\it bad witness}.
The test  picks  a random witness $x$ in 
$W_n$ and evaluates $\rmP_n(x)$.
Two important  characteristics of a pseudo-primality test are 
the run-time {\it complexity} $n\mapsto T(n)$ of the algorithm evaluating
$\mathrm{P}_n$,
 and the {\it density} $n\mapsto \mu(n)$ of bad witnesses. 

To be quite rigorous, we do not need to be able to evaluate
$\mathrm{P}_n$ in deterministic time $T(n)$. 
We are content with a Las Vegas probabilistic  algorithm  that on input
$n$, 
runs in time $T(n)$,  and returns with probability $\geqslant 1/2$
at least one of the 
 following two things
\begin{itemize}
\item[$\bullet$] a proof that  $n$ is composite,
\item[$\bullet$] the value of $\mathrm{P}_n$ at
a random (with uniform probability) element in 
$W_n$.
\end{itemize}

If this is the case, we say that the test $\mathrm{P}$
has complexity $n\mapsto T(n)$ and  density $n\mapsto \mu(n)$.
See   \cite[page 256]{PAPA} 
for the  definition of a Las Vegas algorithm.

\secstar{The Miller-Rabin test} We assume $n$ is odd. The set  $W_n$ of
witnesses for the Miller-Rabin test 
is $(\ZZ/n\ZZ)^*$. The associated map  \[\rmMR_n : (\ZZ/n\ZZ)^*
\rightarrow \{\composite, \mprime\}\] 
is defined by 
$\rmMR_n(x)= \mprime$ if and only if  $x^m=1$ or 
$x^{m2^i}=-1$ for some $0\leqslant i<k$. Here $m$ is
the largest odd divisor of $n-1$ and $n-1=m2^k$.
We call $\rmMR_n$ a {\it Miller-Rabin map}.
It is clear that if $n$ is prime then $\rmMR(x)=\mprime$
for every $x$ in $W_n$.
In case $n$
is composite, the density $\mu_{\rmMR}(n)$  of bad witnesses
is bounded from above
by $1/4$ (see \cite[Theorem $2.1$]{Schoof}).
It will be important for us that  this density 
is actually bounded from above by $2^{1-t}$ (see \cite[proof of Theorem
$2.1$]{Schoof}) 
where $t$ is the number
of prime divisors of $n$. The complexity $T_{\rmMR}(n)$
is bounded from above by $(\log n)^{2+o(1)}$
using fast exponentiation and fast arithmetic.
If we  run $k$  independent Miller-Rabin tests, the probability
of missing a composite number is $\leqslant 4^{-k}$  
and the complexity is $k(\log n)^{2+o(1)}$.

\secstar{A faster pseudo-primality test}\label{sec:gt}
In this article we  prove the following theorem.

\begin{theorem}[A faster test]\label{th:main}
There exist a function $\varepsilon : \RR \rightarrow \RR$
in the class $o(1)$
and a probabilistic algorithm (described in Section~\ref{sec:a-th-algorithm})
that takes as input an odd integer $n$ and an integer
$\lambda$ such that $1\leqslant \lambda\leqslant {\log  n}$, 
runs in time
\[T = (\log n)^{2+\varepsilon(n)} \lambda^{\frac{1}{2}+\varepsilon(\lambda)},\] an returns 
$\mprime$ always if $n$ is prime, and 
with probability \[\leqslant 2^{-\lambda}\] if $n$ is composite.
\end{theorem}

This algorithm  achieves the security of $\lambda/2$ Miller-Rabin tests
at the cost of $\lambda^{1/2+o(1)}$ such tests. The two main ingredients of
our test are the {\it product} of pseudo-primality tests
and a primality criterion involving an  extension of the ring $\ZZ/n\ZZ$.

\secstar{Products}\label{sec:products}
We introduce   the associative  composition law 
\[\vee : \{\composite , \mprime\}\times  \{\composite , \mprime\}
\rightarrow  \{\composite , \mprime\}\] with table

\begin{center}
\begin{tabular}{|c||c|c|}
\hline
$\vee$ &composite &  prime\\
\hline
\hline
composite & composite & composite\\
\hline
prime & composite & prime  \\
\hline
\end{tabular}
\end{center}

\medskip

Let $r\geqslant 2$ be an integer and let
$\mathrm{P}_n^i : W_n^i \rightarrow \{ \composite, \mprime \}$ be $r$
pseudo-primality tests. One defines the product test 
\begin{displaymath}
\rmP_n=\vee_{1\leqslant i\leqslant r} \rmP_n^i
\end{displaymath}
 as 
\begin{displaymath}
\xymatrix@!0 @R=0.6cm @C=5cm{
**[r]\mathrm{P}_n \quad : & W_n= W_n^1\times W_n^2\times \dots \times W_n^r
\ar@{->}[r] & \{\composite, \mprime \} \\
 & (x_1,\dots,x_r)  \ar@{|->}[r] & \vee_{1\leqslant i \leqslant r} \mathrm{P}_n^i(x_i).
}
\end{displaymath}

A witness for $\rmP$ is an $r$-uple of witnesses, one
 for each
of the $r$ tests $\rmP_n ^1$, \ldots, $\rmP_n^r$. For $n$ composite, 
a witness is bad if and only if all its $r$ coordinates are bad 
witnesses. So the density of bad witnesses is the product of
all the densities  for every tests. And the complexity
is bounded 
by the sum of all $r$ complexities, times $\lceil \log_2 r\rceil +1$.
This last factor is natural when chaining Las Vegas algorithms.
In order  to make sure that the resulting algorithm still
succeeds with probability $\geqslant 1/2$ we must repeat a little bit
every step.
As a special case, we consider 
the  $r$-th power $\vee^r \rmP$ of a single test
$\rmP$ with complexity $T$ and density $\mu$.
The density of bad witnesses for
 $\vee^r \rmP$ is  equal to $\mu^r$,
 and its complexity is  $r\times T\times (\lceil \log_2 r\rceil +1)$. 

\secstar{A compositeness criterion}
The test in Theorem~\ref{th:main} is based on the following
compositeness criterion.

\begin{theorem}[Compositeness criterion]\label{th:fermat}
Let $n\geqslant 2$ be an  integer. 
Let $S\supset \ZZ/n\ZZ$ be a faithful, finite,  associative, commutative
$\ZZ/n\ZZ$-algebra  with unit.  
Let  $\sigma$ be an $\ZZ/n\ZZ$-endomorphism of
$S$.
Let $\Omega \subset S$ be a subset of $S$ such that the smallest
$\ZZ/n\ZZ$-subalgebra of $S$  containing $\Omega$
and stable under the action of $\sigma$
is $S$ itself.
Assume $\omega^n=\sigma(\omega)$ for every $\omega$ in $\Omega$.
If $n$ is prime, then for every  $x$ in $S$ we have
$x^n=\sigma(x)$. 
\end{theorem}
\begin{myproof}
Let $T$ be the subset of $S$ consisting of all $x$
such that $x^n=\sigma (x)$. Clearly $T$ contains
$\Omega$. 
If $n$ is prime,  then 
$T$  contains $\ZZ/n\ZZ$ and
is stable under  addition, multiplication, 
and action of $\sigma$. So $T=S$ and we have
$x^n=\sigma(x)$ for every $x$ in $S$. 
\end{myproof}
\medskip

Theorem~\ref{th:fermat} provides a compositeness criterion since
the existence of an $x$ in $S$
such that $x^n\not =\sigma(x)$
implies that $n$ is not a prime.
We call the associated pseudo-primality test a {\it Galois  test}.
The set $W_n$ of witnesses is the group $S^*$ of units in $S$.
The map $\rmP_n$ is defined by $\rmP_n(x)=\mprime$ 
if $\sigma(x)=x^n$ and $\rmP_n(x)=\composite$ otherwise.
In that situation, we call $\rmP_n$ a {\it Galois map}.
In case $n$ is composite, those  $x$ in $S$ for which 
\begin{equation}\label{eq:eq}
x^n=\sigma(x)
\end{equation} are the  {\it bad witnesses}.

\secstar{Plan}
We will show in Section~\ref{sec:ringextensions} that  one
can  bound from above the density of  bad
witnesses among the units of the algebra $S$ in Theorem~\ref{th:fermat}, at least when $S$ is a cyclic
extension of $\ZZ/n\ZZ$. We will use the Galois module
structure of the unit group
of such an  extension. The resulting
pseudo-primality test is presented an analyzed  in Section~\ref{sec:eff}.
Section~\ref{section:construct-algI} explains how to efficiently construct the cyclic $\ZZ/n\ZZ$-algebras required
by our test.  
Theorem~\ref{th:main} is proven
in Section~\ref{sec:a-th-algorithm}.
Implementation details are given in 
Section~\ref{sec:an-algorithm}. We present  the results
of our experiments in Section~\ref{sec:experiments}.

\secstar{Context}
There exist many (families of) algorithms for prime detection.
A recent survey  
can be found in Schoof's article
 \cite{Schoof}. 
The first polynomial time deterministic algorithm for 
distinguishing prime numbers
from composite numbers is due to Agrawal, Kayal and Saxena \cite{aks}.
An improvement of this algorithm,  due to Lenstra and Pomerance \cite{LP},
has deterministic complexity $(\log n)^{6+o(1)}$. This is the best known 
unconditional
result for deterministic algorithms. There exists
a deterministic algorithm with complexity  $(\log n)^{4+o(1)}$
under  the
generalized Riemann hypothesis, as observed by Miller in  \cite{mil}.
Dan Bernstein has found \cite{B} a Las Vegas probabilistic algorithm with
complexity  $(\log n)^{4+o(1)}$. See also Avanzi and Mih\u{a}ilescu
\cite{AvM}.
The correctness and running time of this algorithm does not depend on the truth of any unproved conjecture. It is unconditional.

\secstar{Notation}
In this paper, the
notation $\Theta$ stands for a positive absolute constant.
Any statement containing
this symbol becomes true if the symbol is replaced in every occurrence by some 
large enough real number. Similarly,
the notation $\varepsilon(x)$ stands for a real  function of the real
parameter
$x$ alone, belonging to the class $o(1)$. 

\section{Cyclic  extensions of $\ZZ/n\ZZ$}\label{sec:ringextensions}

Let $n\geqslant 3$ be an odd integer and set
$R=\ZZ/n\ZZ$. 
A {\it cyclic} extension of $R$ is a Galois
extension $S$ of $R$ in the sense of 
\cite[Chapter III]{DI},
with finite cyclic Galois group
$\cG$.  We denote by
$d$ the order of $\cG$, and 
let $\sigma$ be a generator of it.
The Galois property implies 
 \cite[Chapter III, Corollary 1.3]{DI}  that $S$
is a projective $R$-module of constant
rank $d$. Since $R$ is semi-local we deduce 
\cite[II.5.3, Proposition 5]{BourbAC57} 
that $S$ is free of rank $d$.
The  sub-algebra $S^\cG$
consisting of elements
in $S$ fixed by $\sigma$ is $R$ itself \cite[Chapter III,
Proposition 1.2]{DI}. And $S$ is a separable $R$-algebra in the sense
that it is
projective as a module over $S\otimes_RS$. We deduce 
\cite[Theorem 2.5.]{AB}
that $S$
is an unramified extension of $R$. And $S$ is a 
free
$R[\cG]$-module of rank $1$.  Equivalently there
exists a normal basis \cite[Theorem 4.2.]{CHR}.
In this section we study the group of units 
of such an algebra and count the solutions to Equation~(\ref{eq:eq})
in it.
In Paragraph~\ref{sec:GM} we localize at a prime
$p$ and we study the Frobenius
action on the residue algebra. We decompose the unit group as a direct product.
The $p$-part is studied in Paragraph~\ref{sec:UU}, and 
the prime to $p$-part is studied in Paragraph~\ref{sec:V}. In Paragraph~\ref{sec:BW}
we deduce an estimate for the number of bad witnesses.
We refer to the book by DeMeyer and Ingraham \cite{DI} for general properties of Galois extensions,
and to Lenstra  \cite{Lgalois,Lbour} for their use in the context of primality testing.

\subsection{The structure
of $S^*$ as a $\ZZ[\cG]$-module}\label{sec:GM}

 We write  $n=\prod_{p}p^{v_p}$
the prime decomposition of $n$.
If  $p$ and $q$ are two distinct 
 primes dividing $n$, then $p^{v_p}S+ q^{v_q}S=S$. Furthermore,
 the intersection
of all $p^{v_p}S$ for $p$ dividing $n$ is zero. So $S$ is isomorphic
to the product 
\[\prod_{p|n} S/p^{v_p} S= \prod_{p|n} S_p,\]
and this decomposition  is an isomorphism of $\ZZ[\cG]$-modules.
So we can and will assume now that
$n=p^v$ is a prime power.

We set
$\bL=S/pS$ and $\bK=R/pR=\ZZ/p\ZZ$. 
Since $pS\cap R=pR$, the ring  $\bL$ is a faithful
$\bK$-algebra.
The $R$-automorphism 
$\sigma : S\rightarrow S$ induces a $\bK$-automorphism of $\bL$ that we call
$\sigma$ also. 
The $\bK$-algebra   $\bL$ has dimension $d$ and is Galois with group
$\cG$ \cite[Proposition 2.7.]{Lbour}.
From   $\bK = \bL^\cG$ we deduce  
\cite[Chapitre 5, paragraphe 1, num{\'e}ro 9, proposition 22]{BourbAC57} 
that   $\bL$ is integral over  $\bK$. Let  $\pgot$ be a prime ideal in $\bL$.
The intersection  $\pgot\cap \bK$ is a prime ideal in  $\bK$, so it is equal to
 $0$. Since  $0$ is maximal in $\bK$, the ideal  $\pgot$ is 
maximal in  $\bL$ 
\cite[Chapitre 5, paragraphe 2, num{\'e}ro 1, Proposition 1]{BourbAC57}.
Thus  $\bL$ is a ring of  dimension $0$. Since $\bL$ is noetherian,
 it is an artinian ring \cite[Chapitre 4, paragraphe 2, num{\'e}ro 5, Proposition 9]{BourbAC57}. The automorphism
$\sigma$ acts  transitively on the set of prime ideals in $\bL$
\cite[Chapitre 5, paragraphe 2, num{\'e}ro 2, Th{\'e}or{\`e}me 2]{BourbAC57}. 
We denote by  $\cG^Z$ (resp. $\cG^T$) the decomposition  group  (resp. inertia group)
of all these prime ideals.  The Galois property \cite[Proposition 1.2]{DI} implies that the
inertia  group is trivial.
Let  $f$ be the order of $\cG^Z$.
We check that $d=fm$ where  $m$ is the number of   prime ideals in  
$\bL$. Let  $\pgot_0$, $\pgot_1$, \ldots, $\pgot_{m-1}$
be all these prime ideals. 
They are pairwise comaximal: for $i\not =  j$ we have $\pgot_i+\pgot_j=\bL$.
The radical of  $\bL$ is
\[\Ngot = \bigcap_{0\leqslant i\leqslant m-1}\pgot_i =\prod_{0\leqslant i\leqslant m-1}\pgot_i=0,\]
because $\bL$ is unramified over $\bK$.
So the map
\[\bL\longrightarrow  \prod_{0\leqslant i\leqslant m-1}\bL/\pgot_i\]
is an isomorphism of $\ZZ[\cG^Z]$-modules.
For every $i$ in $\{0,1, \ldots, m-1\}$, the  decomposition
group $\cG^Z$ is  isomorphic to the group of  $\bK$-automorphisms
of   the residue field $\bM_i=\bL/\pgot_i$  \cite[Chapitre 5, paragraphe 2, num{\'e}ro 2, Th{\'e}or{\`e}me 2]{BourbAC57}.  
The Frobenius automorphism
$\Phi_i$ of $\bM_i=\bL/\pgot_i$ is the reduction modulo $\pgot_i $ 
of some power $\sigma^{z_im}$ of $\sigma$ generating $\cG^Z$. 
Especially, for every $a$
in $\bL$, one has $\sigma^{z_0m} (a)=a^p\bmod \pgot_0$ for some integer 
$z_0$. We let $\sigma$
act on the above congruence and deduce that $z_0=z_1=\cdots=z_{d-1}
\bmod f$ because $\sigma$ acts transitively on
the set of primes. So there exists a prime to $f$ integer
  $z$ 
 such that for every element $x$ in $\bL$ we have
\begin{displaymath}
  x^p =\sigma^{zm}(x)\,.
\end{displaymath}

We set
\[\UU = \{x\in S | x\equiv 1 \bmod p\}.\]
This is a subgroup of the group  $S^*$ of units in $S$, and even
a  $\ZZ[\cG]$-module. We have an exact sequence
of $\ZZ[\cG]$-modules
\[1\rightarrow \UU\rightarrow S^*\rightarrow (S/pS)^*\rightarrow 1.\]

While the  group $\UU$ is a $p$-group, the group $(S/pS)^*=\bL^*$
has order prime to $p$. So $\UU$ is the $p$-Sylow subgroup
of $S^*$. We denote by $\VV$ the product of all $q$-Sylow subgroups
of $S^*$ for $q\not =p$. Then 
\begin{equation}\label{eq:decomp}
S^*=\UU\times \VV
\end{equation} and this
 decomposition  is an isomorphism of $\ZZ[\cG]$-modules
because both $\UU$ and $\VV$ are characteristic subgroups
of $S^*$. Furthermore,
   $\VV$ is isomorphic to $(S/pS)^*$ as a $\ZZ[\cG]$-module.
We study either factors separately.

\subsection{The structure of $\UU$}\label{sec:UU}

The two maps
\[
\xymatrix@!0 @R=0.6cm @C=1.5cm{
\Log  : & **[l]\UU  \ar@{->}[r] & **[r]p S \\
& **[l]x \ar@{|->}[r] & **[r]\displaystyle \Log (x)=-\sum_{k\geqslant 1}\frac{(1-x)^k}{k}}
\]
and
\[
\xymatrix@!0 @R=0.6cm @C=1.5cm{
\Exp  : & **[l]pS  \ar@{->}[r] & **[r]\UU \\
& **[l]x \ar@{|->}[r] & **[r]\displaystyle \Exp (x)=1+\sum_{k\geqslant 1}\frac{x^k}{k!}}
\]
are well defined. They are indeed polynomial maps (recall
that  $p$ is odd).
In particular, both maps  are equivariant for the action of $\cG$. So 
$\Log$  is an  isomorphism
between the $\ZZ[\cG]$-modules $(\UU,\times)$ and $(pS,+)$. And  $\Exp$  is the 
reciprocal map. 

\subsection{The structure of $\VV$}\label{sec:V}

Let $\pgot$ be a prime in $S$ above $p$. We set 
$\bM=S/\pgot$. Recall that  \[pS=\prod_{0\leqslant k\leqslant m-1}\sigma^{k}(\pgot),\]
and there exists  a prime to $f$ integer
  $z$ 
 such that for every element $x$ in $S$ we have
\begin{displaymath}
  x^p =\sigma^{zm}(x)\,\bmod p.
\end{displaymath}
Let  $1\leqslant t\leqslant f-1$ be the  inverse of $z$ modulo $f$.
Note that if $f=1$, we have $z=t=0$.
We turn $\bM^m$ into a $\ZZ[\cG]$-module by setting 
\begin{equation}\label{eq:action}
\sigma.(x_0, x_1, \ldots, x_{m-1})=(x_1, x_2, \ldots, x_{m-1},x_0^{p^t}).
\end{equation}

The map
\[
\xymatrix@!0 @R=0.6cm @C=1.5cm{
 **[l]S/pS  \ar@{->}[r] & **[r](S/\pgot S)^m \\
 **[l]x \ar@{|->}[r] & **[r]\left(\sigma^k (x) \bmod \pgot\right)_{0\leqslant k\leqslant m-1}}
\]
is an isomorphism of $\ZZ[\cG]$-module between $S/pS$ and 
$\bM^m$. So $\VV$ and $(\bM^*)^m$ are isomorphic as $\ZZ[\cG]$-modules.

\subsection{Counting bad witnesses}\label{sec:BW}

We now show that in many  cases one
can  bound from above the density of bad
witnesses among the units of $S$. 

\begin{theorem}[Density of bad witnesses]\label{th:dens}
Let $A > 2$ and $B\geqslant 3$ be two real numbers.
 Let $n\geqslant 3$ be an    integer. Assume that every 
prime dividing $n$
is bigger than or equal to $B$. %Set $R=\ZZ/n\ZZ$.
Assume that $n$ is not a prime power.
Let $S\supset \ZZ/n\ZZ$ be a cyclic  
$(\ZZ/n\ZZ)$-algebra of dimension
$d$. Let $\sigma$ be a generator of the Galois group $\cG$.
Assume that $n$  has a  prime power
divisor $p^v$  satisfying 
\begin{equation}\label{eq:vpan}
v\log p \geqslant \frac{A\log n}{d}.
\end{equation}
Then the density \[\mu_S=\frac{\#\{x\in S^* | \sigma(x)=x^n\}}{\#S^*}\] of bad witnesses among the units of $S$ is such that

\begin{equation}\label{eq:3}
\mu_S\leqslant {p^{-\frac{vd}{2}(1-\frac{2}{A}-\frac{4}{B})}}\leqslant n^{-\frac{A}{2}(1-\frac{2}{A}-\frac{4}{B})}.
\end{equation}

\end{theorem}

\begin{myproof}
We  count the solutions to Equation~(\ref{eq:eq}) in $S^*$. 
Since $S$ is isomorphic to the product
of all $S_p$ for $p$ a prime dividing  $n$, we fix such 
a prime $p$  and count the solutions to Equation~(\ref{eq:eq})
in $S_p^*$. 
Using the decomposition in Equation~(\ref{eq:decomp})
we then reduce to counting  solutions in the subgroups $\UU$ and $\VV$.

If $x \in \UU$ is a solution to Equation~(\ref{eq:eq}) then
$x^{n^d}=x$. Since $\UU$ is a $p$-group and $p$ divides $n$
we deduce that $x=1$.

According to Section \ref{sec:V}, the $R[\cG]$-module
$\VV$ is isomorphic to $[(S/\pgot S)^*]^m$ where 
$m$ is the number of prime ideals in $S$ above
$p$, and $\pgot$ is one of them, and the action of $\cG$ is given
by Equation~(\ref{eq:action}). It is clear that any solution $x$
to Equation~(\ref{eq:eq}) in the latter $R[\cG]$-module is characterized
by its first coordinate $x_0$ and this coordinate must be a
$|n^m-p^t|$-th root of unity in the field $S/\pgot S$. Since the latter
field has cardinality  $p^f$ we deduce that the number of solutions
to  Equation~(\ref{eq:eq}) in $\VV$ is  
\begin{displaymath}
\gcd(n^m-p^t,p^f-1).
\end{displaymath}

The density of bad witnesses is thus
\begin{equation}\label{eq:2}
\mu_S=\prod_{p|n} \frac{\gcd(n^m-p^t,p^f-1)}{(p^f-1)^mp^{(v-1)d}},
\end{equation}
\noindent where the integers $f,m,v$ and $t$ depend on $p$. This
density
is bounded from above by any term in the product~(\ref{eq:2}).
So let  $p$ 
be a  prime divisor
of $n$ such that 
 $v\log p \geqslant \frac{A\log n}{d}$.
Let $m$ be the number of prime ideals  in $S$
above $p$. 

We first assume that 
 $m\geqslant 2$,
so  $p$ splits in $S$. 
Then the density of bad witnesses is bounded from above by
${1}/{(p^f-1)^{m-1}p^{(v-1)d}}$.  We check that
\begin{equation}\label{eq:beta}
N-1\geqslant N^{(1-\frac{2}{B})},
\end{equation}
for every integer $N\geqslant B$. So $p^f-1 \geqslant p^{{f}(1-\frac{2}{B})}$.
Since  
$m-1\geqslant m/2$, we find 
\[\mu_S \leqslant {1}/{p^{\frac{d}{2}(1-\frac{2}{B})+(v-1)d}}.\]
The result follows.

We now assume  that $m=1$,  so $p$ is inert in $S$ and  $f=d$.
We first prove  the following inequality
\begin{equation}\label{eq:5}
\gcd(n-p^t,p^d-1) \leqslant {np^\frac{d}{2}}.
\end{equation}
Indeed, if $1\leqslant t\leqslant\frac{d}{2}$, Inequality~(\ref{eq:5}) is granted
 because
$1\leqslant |n-p^t|\leqslant \max (n,p^t)\leqslant np^t$.
In case  $\frac{d}{2} <  t \leqslant d-1$, we call  $w$ the unique 
integer in $[1,d[$ that is congruent to  $-t$
modulo $d$. We have
\begin{equation}\label{eq:4}
\gcd(n-p^t,p^d-1)=\gcd(np^w-1,p^d-1).\\
\end{equation}
 Since $w \leqslant  (d-1)/2$, the right hand side of (\ref{eq:4}) is bounded
 from above by $np^{\frac{d}{2}}$ as was to be shown. 
So Inequality~(\ref{eq:5})
holds true  in either  case, and Inequality~(\ref{eq:3}) follows using
Equation~(\ref{eq:2}), Equation~(\ref{eq:vpan}),  and Inequality~(\ref{eq:beta}). 
\end{myproof}

\section{An efficient pseudo-primality test}\label{sec:eff}

A consequence of Theorem~\ref{th:dens} is
that a  compositeness criterion as   Theorem~\ref{th:fermat}, when
implemented with
a cyclic $(\ZZ/n\ZZ)$-algebra of dimension $d$,
is efficient, provided   $n$ has a large prime power divisor
$p^v$.
On the other hand, we saw in Section~\ref{sec:intro} that
the  Miller-Rabin test is efficient when $n$ has many prime
divisors. 
Combining these
two tests  we can construct a new probabilistic pseudo-primality test
that takes advantage of either situation.

Fix  two real numbers $A$ and $B$ such that $A>2$ and 
$B\geqslant 4A/(A-2)$. In particular $B>4$. Set $C=1-2/A-4/B$ and
note that $C$ is  positive.

Let  $n$ be a positive  integer. We assume $n$
is not a prime power, and every prime dividing $n$
is bigger than or equal to $B$.
 We choose
 two positive    integers $r$ and 
$d$ and we construct a pseudo-primality test which is the product of
$r$ Miller-Rabin tests and a Galois test of dimension
$d$.
We let $\delta = \log (d/A)/\log \log n$ so \[{d}=A(\log n)^\delta.\]
We let $\rho = \log (2A^{-1}r\log 2)/(\log \log n)$ so 
\[{r}=\frac{A(\log n)^\rho}{2\log 2}.\]
We assume 
\begin{equation}\label{eq:co1}
\left(1-\frac{A}{d}\right)(\log n)^{\delta+\rho}\leqslant C\log n,
\end{equation}
or equivalently
\begin{displaymath}
dr\left(1-\frac{A}{d}\right)\leqslant \frac{A^2C\log n}{2\log 2}.
\end{displaymath}

We call $\rmP_1 : ((\ZZ/n\ZZ)^*)^r \rightarrow \{ \composite, \mprime \}$
 the product  of $r$ Miller-Rabin maps.
And $\rmP_2 : S^* \rightarrow \{ \composite, \mprime \}$  a 
Galois map as in Theorem~\ref{th:fermat},
associated with a
cyclic  algebra of dimension $d$.
We set $\rmP=\rmP_1\vee \rmP_2$.
The density  of bad witnesses for $\rmP$ is bounded from above
by the densities of bad witnesses for $\rmP_1$ and $\rmP_2$.
Let $p^v$ be the largest prime power dividing $n$.
We set $\pi = \log (v\log p)/\log \log (n)$, so \[\log p^v =(\log n)^\pi.\]
The number $t$ of prime divisors of $n$ satisfies 
\[t > (\log n)/(v\log p)=(\log n)^{1-\pi}.\]
If \[\delta+\pi \geqslant 1,\] then
$v\log p \geqslant \frac{A\log n}{d}$, and, according to
Theorem~\ref{th:dens}, the density  of bad
witnesses for $\rmP_2$ is bounded from above by 
\begin{equation}\label{eq:b1}
p^{-\frac{vd}{2}(1-\frac{2}{A}-\frac{4}{B})}=\exp(-\frac{A}{2}(1-\frac{2}{A}-\frac{4}{B})(\log n)^{\delta+\pi}).
\end{equation}
On the other hand,  the density of bad witnesses for every
Miller-Rabin test is $\leqslant {2^{-t+1}}$.
The density of bad witnesses for $r$ such tests is at most
\begin{equation}\label{eq:b2}
2^{-{r(t-1)}}\leqslant \exp(-\frac{A}{2}(1-\frac{1}{t})(\log n)^{1+\rho-\pi}).
\end{equation}
  
Although we do not know the value of $\pi$, we can deduce
from Equations~(\ref{eq:b1}) and~(\ref{eq:b2}) an upper bound
for the density  of bad witnesses of the product test
$\rmP=\rmP_1\vee \rmP_2$.

If $\pi$ lies in $[0,1-\delta[$ then Equation~(\ref{eq:b1}) gives nothing and Equation~(\ref{eq:b2})
gives an upper bound
\begin{displaymath}
\exp(-\frac{A}{2}(1-\frac{A}{d})(\log n)^{\rho+\delta}),
\end{displaymath}
 for the density of bad witnesses for $\rmP_1$.

If $\pi$ lies in $[1-\delta,1]$ then Equation~(\ref{eq:b1}) gives  an upper bound 
\begin{displaymath}
\exp(-\frac{A}{2}(1-\frac{2}{A}-\frac{4}{B})\log n),
\end{displaymath}
 for the density of bad witnesses
for $\rmP_2$. Using Inequality~(\ref{eq:co1}) we find the upper bound
\begin{displaymath}
\exp(-\frac{A}{2}(1-\frac{A}{d})(\log n)^{\rho+\delta}),
\end{displaymath}
in that case.

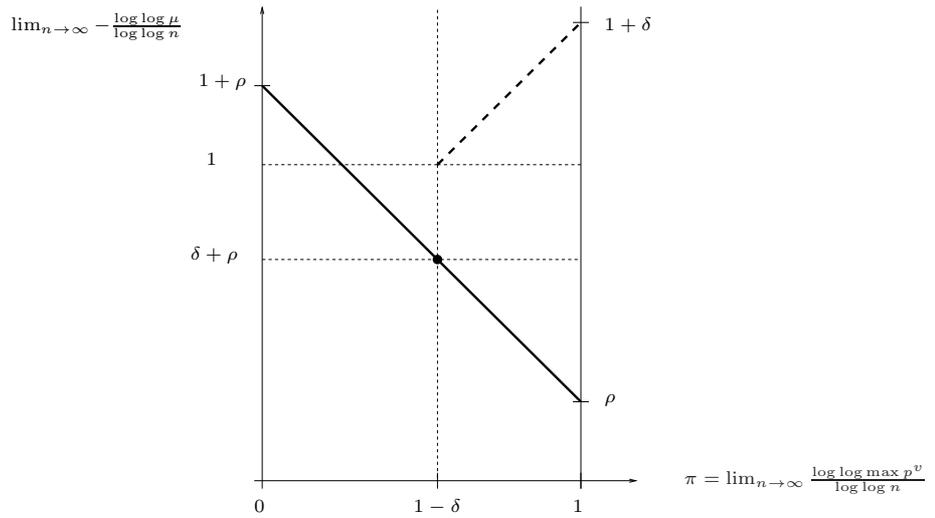
\begin{figure}
  \begin{center}
    \input{courbes.pspdftex}
  \end{center}
  \caption{The Miller-Rabin (continuous)  and Galois (dashed) densities.}
  \label{fig:qc}
\end{figure}

This discussion is illustrated in Figure~\ref{fig:qc} where the
continuous line 
is
the exponent of $\log n$ in Equation~(\ref{eq:b2}),
the dashed 
line  is the exponent of $\log n$ in
Equation~(\ref{eq:b1}),
and the bullet  is the minimum of the maximum of the two functions.

\begin{theorem}[Density of the composed test]\label{th:efftest}
Let  $A$ and $B$ be two real numbers such that $A>2$ and 
$B\geqslant 4A/(A-2)$.  Let 
\begin{equation}\label{eq:C}
C=1-2/A-4/B.
\end{equation}
Let $n$ be an    integer
that is not a prime power. Assume that $n$ has no prime divisor 
smaller than $B$.
Let  $r$ and 
$d$ be  two positive    integers such that 
\begin{equation}\label{eq:t}
dr\left(1-\frac{A}{d}\right)\leqslant \frac{A^2C\log n}{2\log 2}
\end{equation}
and let
$\rmP$ be the composite test of 
$r$ Miller-Rabin tests and one  Galois test of dimension
$d$.
The density of bad witnesses for $\rmP$ is bounded from above by 
 \[\leqslant 2^{-\frac{rd}{A}(1-\frac{A}{d})}.\]
\end{theorem}

Taking $A=2.1$, $B=1000$, and $d\geqslant 16$, we have
$C\geqslant 0.043619$ and 
we obtain a density
$\leqslant 2^{-0.41369 rd}$ provided
$rd\leqslant 0.13875 \log n$.

Taking $A=4$, $B=1000$, and $d\geqslant  16$, we have
$C\geqslant 0.496$ and 
we obtain a density
$\leqslant 2^{-0.18 rd}$ provided
$rd\leqslant 5.72 \log n$.

We note that the complexity of such  a composed
test is $(\log n)^{2+\varepsilon (n)}(r+d^{1+\varepsilon (d)})$ under the condition
that arithmetic operations in the $\ZZ/n\ZZ$-algebra $S$ can be performed in
quasi-linear time in the degree $d$.
It is asymptotically
optimal to take $d$ and $r$ as close as possible. 
We thus  prove Theorem~\ref{th:main} provided we can 
efficiently construct a Galois extension of $\ZZ/n\ZZ$ with degree $d$ in
some  interval $[k,k^{1+\varepsilon (k)}]$.
This is the purpose of the next Section~\ref{section:construct-algI}.

\secstar{Heuristics}
There are many possible choices for the parameters $A$, $B$, $r$ and  $d$ when
using Theorem~\ref{th:efftest}. We will explain in
Section~\ref{sec:an-algorithm} how to choose them optimally. Here we
just
collect a few simple minded observations on what could be a reasonable choice.
We  take
\begin{equation}\label{eq:cB}
B=8000. 
\end{equation}
Taking a too large $A$ is pointless. We recommend 
\begin{equation}\label{eq:cA}
2<A\leqslant 48.
\end{equation}
In case we have a bigger value of $A$ it will be more efficient to
take smaller values for  $r$ and $d$
and repeat the whole test.
We  also suggest  that 
\begin{equation}\label{eq:cdA}
d\geqslant 2A,
\end{equation}
otherwise we would better use $r$ Miller-Rabin tests only, and obtain better security at  lower cost.
It is reasonable also to have
\begin{equation}\label{eq:cdr}
d\leqslant r,
\end{equation}
because the $r$ Miller-Rabin  tests and the one Galois test have similar effect on the security.
So the time devoted to  the $r$ Miller-Rabin  tests  should not be smaller than the time
devoted to the Galois test.
Assume we want to bound from above the error probability by $2^{-\lambda}$ for some integer $\lambda$.
We must have
\begin{equation}\label{eq:cL1}
\lambda \leqslant \frac{rd}{A}(1-\frac{A}{d}).
\end{equation}
And  we should have 
\begin{equation}\label{eq:cL2}
\frac{rd}{A}(1-\frac{A}{d})\leqslant 2\lambda,
\end{equation}
in order not to waste time.

We deduce from Equations~(\ref{eq:cdr}),  (\ref{eq:cL2}),
(\ref{eq:cdA}),
and~(\ref{eq:cA}) that 
\begin{equation}\label{eq:cd}
d\leqslant 2\sqrt{A\lambda}\leqslant 14\sqrt\lambda.
\end{equation}

We deduce from Equations~(\ref{eq:cL1}), (\ref{eq:t}), (\ref{eq:C}),
and~(\ref{eq:cA}),  that 
\begin{equation}\label{eq:cL3}
\lambda\leqslant (0.9995A-2)\frac{\log_2 n}{2}\leqslant 23\log_2 n.
\end{equation}

Under the reasonable hypotheses above, 
the smallest possible value for $A$ when
 applying Theorem~\ref{th:efftest} is  thus
\begin{displaymath}
\left(2+\frac{2\lambda}{\log_2(n)}\right)/0.9995.
\end{displaymath}
So we recommend to take 
\begin{equation}\label{eq:recA}
A=\left(2+\frac{2\lambda}{b-1}\right)/0.9995,
\end{equation}
 where 
 \[b = \lfloor \log_2(n) \rfloor+1,\]
is the
number of bits of $n$.
\section{Constructing algebras}\label{section:construct-algI}

In this section we prove the following theorem.

\begin{theorem}[Constructing algebras]\label{th:const}
There exist a function $\varepsilon : \RR\rightarrow
\RR$ in the class $o(1)$ and a probabilistic (Las Vegas)
algorithm
that takes as input  an odd  integer $n$ and an integer 
$k$ such that $1\leqslant k\leqslant \log n$, 
runs in time $(\log n)^{2+\varepsilon (n)}$, and returns  with probability $\geqslant 1/2$ at least one of the
 following two data
\begin{itemize}
\item[$\bullet$] A proof that  $n$ is composite,
\item[$\bullet$] A cyclic algebra $S$
over $\ZZ/n\ZZ$ with degree
$d$ and Galois group $\cG = \langle \sigma\rangle $
such that 
\begin{equation}\label{eq:borned}
k\leqslant d\leqslant k^{1 +\varepsilon(k)},
\end{equation}
and there exists a basis $\Omega$ of the $\ZZ/n\ZZ$-module
$S$ such that $\sigma(\omega)=\omega^n$ for every
$\omega$ in $\Omega$.
\end{itemize}
Arithmetic operations in $S$ are then performed in deterministic time
$(\log n)^{1+\varepsilon(n)}d^{1+\varepsilon(d)}$.
\end{theorem}

From Theorem~\ref{th:const} and Theorem~\ref{th:efftest}
one can easily deduce Theorem~\ref{th:main}.
We prove Theorem~\ref{th:const} in two steps. 
We first apply a single  Miller-Rabin test to $n$. If $n$ is composite we
shall thus 
detect it  with probability $\geqslant 1/2$  in probabilistic time $(\log
n)^{2+\varepsilon (n)}$. So this copes with the case when $n$ is composite.
We then try to construct an
$(\ZZ/n\ZZ)$-algebra  $S$. 
For the  complexity
analysis of this second step, we can assume that $n$ is prime.

We shall use Kummer theory to construct an extension
of $\ZZ/n\ZZ$ with appropriate degree. This is a classical construction in this context. It appears in \cite{apr,LP} and 
even more  explicitly    in \cite{B,KU}. We first
construct a small cyclotomic extension $\Rc$, then a Kummer
extension $S$ of $\Rc$.
We let $\dc$ be the smallest positive integer such that
the product  $Q$ of all prime integers $q$ such that
$q-1|\dc$ exceeds $k$. According to \cite[Theorem 3]{apr}
we have 
\begin{displaymath}
\dc \leqslant  (\log k)^{\Theta \log\log\log \Theta k}.
\end{displaymath}

We call   $\dk$ the smallest divisor of $Q$ that exceeds
 $k$. We set $d=\dk\dc$.
It is clear that  
$d$  satisfies Inequality~(\ref{eq:borned}).
We first use the algorithms  in \cite{Shoup} to find a degree $\dc$
unitary polynomial $F(X)$ in $\ZZ/n\ZZ[X]$ that is irreducible if $n$ is
prime.
This takes probabilistic time $\dc^{2+\varepsilon(\dc)}(\log n)^{2+\varepsilon(n)}$ that
is $(\log n)^{2+\varepsilon(n)}$.
We set \[\Rc = (\ZZ/n\ZZ)[X]/F(X).\]
We set $x=X\bmod F(X)$ and  call $\sigmac : \Rc \rightarrow \Rc$ the
$(\ZZ/n\ZZ)$-linear map
that sends $x^i$ to $x^{ni}$ for $0\leqslant i\leqslant \dc -1$.
We check that $\sigmac$ is a morphism of $(\ZZ/n\ZZ)$-algebras. This boils down
to checking that $\sigmac(x^i)=x^{ni}$ for $\dc \leqslant i\leqslant
2\dc -2$. This takes time $(\log n)^{2+\varepsilon(n)}$.
It is a matter of linear algebra to  check that the fixed subalgebra 
by $\sigmac$ is $\ZZ/n\ZZ$. It takes time $(\dc)^3(\log n)^{1+\varepsilon(n)}=(\log
n)^{1+\varepsilon(n)}$.
We pick a random $\xc$ in $\Rc$ and check that 
\begin{equation}\label{eq:nrc}
\sigmac ^i (\xc)-\xc \in \Rcs
\end{equation}
for every $0 <  i <  \dc$. If $n$ is prime
then the density of such elements  in $\Rc$ is at least $1/2$.
So finding one of them takes probabilistic time $(\log n)^{2+\varepsilon (n)}$.

We check that $\dk$ divides $n^\dc-1$. We check that $\sigmac^\dc (x)=x$.

We  look for an element $a$ in $\Rcs$
such that
$\zeta = a^{\frac{n^{\dc} -1}{\dk}}$
has exact order $\dk$.
If $n$ is prime, the density of such elements $a$
 in $\Rcs$ is
$\geqslant (\log \log \log n)^{-\Theta}$.
We check that $\sigmac (a)=a^n$.

We  set
 \[S=\Rc[Y]/(Y^{\dk}-a),\]
and  $y=Y\bmod Y^{\dk}-a$. Let $\tau : S\rightarrow  S$ be the 
unique endomorphism of $\Rc$-algebra such that $\tau (y)=\zeta y$.
The fixed subalgebra by $\tau$ in $S$ is $\Rc$.

There exists a unique endomorphism of  $(\ZZ/n\ZZ)$-algebra
$\sigma : S \rightarrow S$ such that 
$\sigma (y)=y^n$ and the restriction of
$\sigma$ to $\Rc$ is $\sigmac$. It is clear that 
$\sigma^{\dc}$ is $\tau$. Restriction to $\Rc$ gives
an exact sequence 
\[1\rightarrow \langle\tau\rangle\rightarrow \langle
\sigma\rangle \rightarrow
\langle \sigmac\rangle \rightarrow 1.\]
So the order of $\sigma$ is $d=\dk\dc$.
Every element in $S$ fixed by $\sigma$ is also fixed
by $\tau = \sigma^{\dk}$. So it belongs to $\Rc$.
But elements in $\Rc$ fixed by $\sigmac$ actually lye in
$\ZZ/n\ZZ$. So 
\begin{equation}\label{eq:fix}
S^\cG=\ZZ/n\ZZ,
\end{equation}
where $\cG$ is the group generated by $\sigma$.
Furthermore,  for every $0<i<\dk$ 
\begin{equation}\label{eq:nrk}
\tau^i(y) -y=(\zeta^i-1)y\in S^*.
\end{equation}

From~(\ref{eq:fix}),  (\ref{eq:nrc}), (\ref{eq:nrk}) and 
\cite[Proposition 1.2]{DI} we deduce that $S$ is a Galois extension
of $\ZZ/n\ZZ$ with group $\cG$.
As for the basis $\Omega$ we can  take the $x^iy^j$
for $0\leqslant i<\dc$ and $0\leqslant j<\dk$.

\secstar{Remark} We  expect \cite[Remark 6.3]{apr} that 
\begin{displaymath}
\dc \leqslant  (2\log \dk)^{1.5 \log\log\log \dk},
\end{displaymath}
for large enough $k$.  This and Equations~(\ref{eq:cd}),
(\ref{eq:cL3})  implies
\begin{equation}\label{eq:bdc}
\dc \leqslant  \left(9+\log b \right)^{1.5\times  \max(1,\log\log\log
  68\sqrt{\log_2 n})},
\end{equation}
where $b$ is the number of bits of $n$.
We shall use this estimate in Section \ref{sec:an-algorithm}.

\section{An algorithm}

It is now possible to specify an algorithm.

\subsection{A theoretical algorithm}\label{sec:a-th-algorithm}

We prove Theorem~\ref{th:main} by describing the algorithm.  The input
consists of a large enough integer $n$ and a bound $\lambda$ such that $1\leqslant
\lambda\leqslant \log n$.  The algorithm
 outputs either that $n$ is composite or that $n$ is a
probable prime.  The probability of missing a composite is at most
$2^{-\lambda}$.

The algorithm is the following.
\begin{enumerate}[label=\roman{*}),topsep=2pt]
  \item Check that $n$ has no prime factor smaller than  $1000$.
  \item Check that $n$ is not a prime power.
  \item Set $k=\max(16,\lfloor \sqrt{\lambda}\rfloor)$
and use the algorithm in the proof of Theorem~\ref{th:const}
to construct a $(\ZZ/n\ZZ)$-algebra $S$ with degree $d$ such that 
$k\leqslant d\leqslant k^{1+\epsilon(k)}$.
  \item Set $r=\lceil {\lambda}/({0.18\times d})\rceil$.
  \item Perform $r$ Miller-Rabin tests. If one of them
fails output $\composite$.
  \item Choose at random a non-zero  $z$ in $S$ and check that
    it is invertible. If it is not, output $\composite$.
  \item Check that $\sigma(z) = z^n$ and output
$\composite$ or $\mprime$ accordingly.
\end{enumerate}

Applying Theorem~\ref{th:efftest} with $A=4$ and $B=1000$
we see that, for large enough $n$,
the algorithm returns $\mprime$  with probability
$\leqslant 2^{-\lambda}$ when $n$ is composite. It
runs in time
$(\log n)^{2+\epsilon(n)}\lambda^{\frac{1}{2}+\epsilon(\lambda)}$
because both $d$ and $r$ are $\leqslant \lambda^{\frac{1}{2}+\epsilon(\lambda)}$.

\subsection{A practical algorithm}\label{sec:an-algorithm}

We let $b$ be the number of bits of $n$.  We assume $\lambda \leqslant
23\log_2 n$ according to Equation~(\ref{eq:cL3}). For higher security we may
just repeat the test.  We set $B=8000$ and
$A=\left(2+\frac{2\lambda}{b-1}\right)/0.9995$ following
Equations~(\ref{eq:cB}) and~(\ref{eq:recA}).

The algorithm of Section~\ref{sec:a-th-algorithm} can be reformulated as
follows.
\begin{itemize}[topsep=2pt]
\item Preliminaries.
  \begin{enumerate}[label=\arabic{*}),start=1,topsep=0pt]

  \item Check that $n$ has no prime factor smaller than  $B$.
  \item Check that $n$ is not a prime power.
  \item Determine the integers $\dc$, $\dk$ and $r$.
  \end{enumerate}\medskip
\item Miller-Rabin tests.
  \begin{enumerate}[label=\arabic{*}),start=4]
  \item Perform $r$ Miller-Rabin tests.
  \end{enumerate}\medskip
\item Construction of the algebra $\Rc$.
  \begin{enumerate}[label=\arabic{*}),start=5]
\item Find an ``irreducible'' polynomial $F(X)$ of degree $\dc$ modulo
  $n$ and construct the algebra $\Rc$.
\item Compute the action of the automorphism $\sigmac$ on every
  $X^i\bmod F(X)$ for $i=0,\ldots,2\dc-2$.
\item Check that the fixed submodule by $\sigmac$ in $\Rc$ is  $\ZZ/n\ZZ$.
\item Find a $u$ in  $\Rc$  such that  $\sigmac^i(u) -u$ is a unit for every $1\leqslant i\leqslant \dc -1$.
  \end{enumerate}\medskip
\item Construction of the algebra $S$.
  \begin{enumerate}[label=\arabic{*}),start=9]
  \item Find an element $a$ in $\Rc$ such that
    $\zeta=a^{\frac{n^\dc-1}{\dk}}$ has exact  order $\dk$. Check that
$\sigmac(a)=a^n$.
  \end{enumerate}\medskip
\item The Galois test.
  \begin{enumerate}[label=\arabic{*}),start=10]
  \item Choose at random a non-zero  $z$ in $S$ and check that
    it is invertible.
  \item Check that $\sigma(z) = z^n$.
  \end{enumerate}
\end{itemize}

We now comment on each of  these steps. 

\subsubsection{Preliminary steps}
\label{sec:preliminaries}

\subsecstar{Step 1: Check that $n$ has no prime factor smaller than $B$}
\label{sec:trial-division-up}
Recall that   $B=8000$.
We  compute once and for all the product of all the primes smaller than $B$
and check that the $\gcd$ with $n$ is equal to $1$. If this is not the case, we
stop and output that $n$ is composite.

\subsecstar{Step 2: Check that $n$ is not a prime power}
\label{sec:check-that-n}
For each integer $d$ between $2$ and $b$, we compute some integer
approximation $\eta$ of the positive real  $\sqrt[d]{n}$  such that 
$|\eta-\sqrt[d]{n}|\leqslant 0.6$ (there exist fast
methods based on Newton iterations for this task). Then we check that
$\eta^d$ is not equal to $n$. Otherwise  we stop and output that $n$ is composite.

\subsecstar{Step 3: Determine the integers $\dc$, $\dk$ and $r$}
\label{sec:determ-best-poss}
We consider all the small integers $\dc$, starting from $1$ and ending at
$\lfloor 
\left(9+\log b\right)^{1.5\times  \max(1,\log\log\log
  68\sqrt{\log_2 n})}\rfloor$ according to Equation~(\ref{eq:bdc}).
For each $\dc$, we enumerate  the divisors
$\dk$ of $n^{\dc}-1$ upper  bounded by 
$\lfloor {2 \sqrt{A \lambda}}/{\dc} \rfloor$
according to Equation~(\ref{eq:cd}).
We set  $d=\dc\times\dk$  and 
 $r = \lceil \lambda\,A / (d-A) \rceil$\,. 

This exhaustive search produces many  $3$-uples ($\dc$, $\dk$,  $r$).
Among these we select  the one with
the smallest estimated cost. The cost estimates are obtained from 
some systematic experiments with  the available
computer  arithmetic (see
Section~\ref{sec:experiments} for our choices in a \textsc{magma}
implementation).

We compare then with the estimated cost of $\lambda/2$ classical Miller-Rabin
tests. If the latter are cheaper, we switch to these classical tests and
output the result, otherwise we go to Step~4.

\subsubsection{Miller-Rabin tests}
\label{sec:miller-rabin-tests}

\subsecstar{Step 4: Perform $r$ Miller-Rabin tests}
\label{sec:perform-r-miller}
Each of these $r$ tests is a classical Miller-Rabin test as described
in Section~\ref{sec:intro}.

\subsubsection{Construction of the algebra $\Rc$}
\label{sec:constr-algebra-rc}
\smallskip

We skip the next four steps when $\dc=1$.

\subsecstar{Step 5: Find a unitary  ``irreducible'' polynomial $F(X)$ of degree $\dc$ modulo $n$}
\label{sec:find-an-irreducible}
We use any efficient probabilistic algorithm  $\cA$
that produces a degree $\dc$ unitary irreducible
polynomial,  with probability $\geqslant 1/2$, provided $n$ is prime. 
For $n$ prime, $\cA$ fails with probability $\leqslant 1/2$. In that case
it returns   nothing.
If $n$ is not prime, then $\cA$ may return either nothing or a unitary 
polynomial of degree $\dc$ in $(\ZZ/n\ZZ)[X]$.

We call $\cB$ the algorithm consisting of $\cA$ followed by a
Miller-Rabin
test.
It returns with probability $\geqslant 1/2$ either a proof that $n$
is not prime or a polynomial of degree $\dc$ in $(\ZZ/n\ZZ)[X]$.
We iterate $\cB$ until we get such an output.

Step 5  thus  provides either a proof of compositeness or
a polynomial which we know to be irreducible in case $n$ is a prime.
As for the choice of $\cA$ we distinguish several cases, for efficiency purposes.
\begin{itemize}
\item When $\dc=2$, we look for an element $o$ with Jacobi Symbol
  $\left(\frac{o}{n}\right)$ equal to $-1$ and we set $F(X)=X^2-o$. Note that
  $o$ is a quadratic non-residue when $n$ is a prime.
\item When $\dc$ divides $n-1$,  we look  for an element $o$
  such that $o^\frac{(n-1)}{\dc}$ has order $\dc$, and we set $F(X)=X^{\dc}-o$.
\item Otherwise, we test   random unitary polynomials $F(X)$  
and we use  the extended Euclidean algorithm to  check that  the
ideal 
$(X^{n^i}-X,F(X))$ in $(\ZZ/n\ZZ)[X]$ is one
for all $i$ from $1$ to $\lfloor \dc/2
  \rfloor$. If we test more than $\log(1/2)/\log(1-1/2d)$ polynomials
$F(X)$, then the probability of success is $\geqslant 1/2$ provided 
$n$ is prime.
\end{itemize}

One may wonder why we incorporate a Miller-Rabin test in the
loop. This is just to guarantee that we leave the loop in due time,
even if $n$ is composite. 
A similar caution should be taken in every  loop occurring in the next
steps. We only detail  this here. In practice these Miller-Rabin test are
completely
useless. Indeed  $n$ is almost known to be prime and there is no risk
that we keep blocked in such a loop.

\subsecstar{Step 6: Compute the action of the automorphism $\sigmac$}
\label{sec:comp-acti-autom}
We set $x=X\bmod F(X)$ and write $x^{i\,n}$ in the polynomial basis $(x^k)_k$,
for $i$ from $0$ to $\dc-1$. This yields a $\dc\times\dc$ matrix over
$\ZZ/n\ZZ$, that we denote $M_{\sigmac}$.  Using this matrix, we can check
that $\sigmac(x^i)=x^{in}$ for $i$ from $\dc$ to $2\,\dc-2$, and $\sigmac^\dc
(x)=x$. If this is not the case, we stop and output that $n$ is composite.

\subsecstar{Step 7: Check that $\sigmac$ fixes $\ZZ/n\ZZ$}
\label{sec:check-that-sigmac}
We try to compute the kernel of $M_{\sigmac}-\Id$, using Gauss
elimination.
It produces either the expected kernel or a zero divisor in
$\ZZ/n\ZZ$.
In the latter case we stop and
output that $n$ is composite. 
Once computed the kernel, we  check that it is equal to $\ZZ/n\ZZ$. 
If it is not the case, we stop and
output that $n$ is composite. 

\subsecstar{Step 8: Find a $u$ in  $\Rc$  such that  $\sigmac^i(u) -u$ is a unit for every $1\leqslant i\leqslant  \dc -1$}
If $n$ is prime then at least half of the elements in $\Rc$ satisfy
the condition.
So we pick at random  $u$ in $\Rc$ and test the condition. We iterate if
it fails. We again add a Miller-Rabin test in the loop to make sure
that it stops with probability $\geqslant 1/2$ even when $n$ is composite.

To check that a non-zero element $z$ in $\Rc$ is  a unit we try to compute an
inverse  using extended Euclidean algorithm. 
If it returns an element $z'$, we just need to
check that $z\,z'=1$. It it fails we know that $n$ is not a prime and we
stop.

\subsubsection{Construction of the algebra $S$}
\label{sec:constr-algebra-s}

\subsecstar{Step 9: Find an element $\zeta$ of exact order $\dk$ in $\Rc$}
\label{sec:find-an-element}
We pick a random $a$ in the algebra $\Rc$ and compute $\zeta =
a^{(n^{\dc}-1)/\dk}$. If $n$ is prime then the density of $a$
such that the corresponding $\zeta$ has exact order $\dk$ is
$\geqslant (\log\log\log n)^{-\Theta}$. 
The test consists of checking that $\zeta^{\dk/q}-1$
is a unit,  for every prime divisor $q$ of $\dk$. We proceed
as in Step 8.

As above, we  add a Miller-Rabin test in the loop to make sure
that it stops with probability $\geqslant 1/2$  when $n$ is composite.

We check that $\sigmac(a)=a^n$ using 
the matrix $M_{\sigmac}$. If this is not the case, 
we know that $n$ is not a prime and we
stop.

\subsubsection{The Galois test}
\label{sec:galois-test}

\subsecstar{Step 10: Choose at random an invertible element in $S$ }
\label{sec:choose-at-random}
We pick a random non-zero
$z$ in $S$ and try to compute the inverse $z'$ of $z$ with
the extended $\gcd$ algorithm.  If the extended $\gcd$ algorithm fails, or $z'
\times z$ is not equal to 1, then
we know that $n$ is not a prime and we
can
stop.

\subsecstar{Step 11: Check that $\sigma(z) = z^n$}
\label{sec:check-that-sigmax}
On the first hand, we compute $z^n$ in $S$ using  fast exponentiation.
On the other hand, we write
\begin{math}
  z = \sum_i z_i\, y^i\,
\end{math}
where $z_i\in\Rc$ and $y=Y\bmod Y^\dk-a$. Then, we compute $\sigma(z)$ as
\begin{displaymath}
   \sum_i \sigmac(z_i)\times y^{in}
\end{displaymath}
where $\sigmac(z_i)$ is computed using the
matrix $M_{\sigmac}$.
Note that $y^{in}$ can be efficiently computed as 
$a^{\alpha} y^{\beta}$ where
$\alpha$ (resp. $\beta$) is the quotient (resp. the remainder)
in the Euclidean division of $in$ by $\dk$.

If  $\sigma(z)$ is not equal to $z^n$, we output that $n$ is
composite. Otherwise, we output that $n$ is a Galois pseudo-prime.

\section{Experiments}\label{sec:experiments}

We first have  determined power functions that best approximate
the
sub-quadratic timings that we have measured for elementary arithmetic
polynomial operations in \textsc{magma v2.18-2}.  In our testing ranges,
\textit{i.e.} $b$ between $512$ and $8192$ bits, $\dc$ between $1$ and $16$
and $\dk$ between $8$ and $1000$, we have obtained  the following upper
bounds for the heaviest steps in the algorithm.
\begin{itemize}
\item Step 4. Computing $r$ Miller-Rabin tests:
  \begin{displaymath}
    T_{\mathrm{MR}}(b, r)=F\times r \times b^{2.6}\,.
\end{displaymath}
\item Step 5. Constructing an ``irreducible'' polynomial of degree $\dc$
  modulo $n$ (worst case):
  \begin{displaymath}
    T_{\mathrm{F}}(b, \dc) = \left\{
      \begin{array}{ll}
        0 & \text{if } \dc = 1\,,\\ 
        F\times \log_2 b \times b^{2.6} & \text{if } \dc = 2\,,\\
        18\,F\times \log_2 \dc \times  \dc^{\,2.2} \times  b^{\,2.4}& \text{for larger } \dc\,.
      \end{array}
    \right.
  \end{displaymath}
\item Step 9. Finding an element $\zeta$ of order $\dk$ in $\Rc$ (worst case):
  \begin{displaymath}
    T_{\mathrm{\zeta}}(b, \dc) = \left\{
      \begin{array}{ll}
        19\,F \times  b^{\,2.4} & \text{if } \dc = 1\,,\\ 
        36\,F \times \dc^{\,2.2} \times b^{\,2.4} & \text{otherwise.}
      \end{array}
    \right.
  \end{displaymath}
\item Step 11. Computing $\sigma(x)$ in $S$:
  \begin{displaymath}
    T_{\sigma}(b, \dc, \dk) = \left\{
      \begin{array}{ll}
        F \times \dk \times b^{\,2.6} & \text{if } \dc = 1\,,\\ 
        10\,F \times (\dc\times\dk) \times b^{\,2.4} & \text{otherwise.}
      \end{array}
    \right.
  \end{displaymath}
\item Step 11 bis. Computing $x^n$ in $S$:
  \begin{displaymath}
    T_{\mathrm{power}}(b, \dc, \dk) = \left\{
      \begin{array}{ll}
        19\,F \times \dk^{\,1.2} \times b^{\,2.4} & \text{if } \dc = 1\,,\\ 
        36\,F \times (\dc\times\dk)^{\,1.2} \times b^{\,2.4} & \text{otherwise.}
      \end{array}
    \right.
  \end{displaymath}
\end{itemize}
For the sake of completeness, we found that the constant $F$ is equal to
$30\times 10^{-{9}}$ seconds on our laptop (based on a \textsc{Intel Core i7
  M620 2.67GHz} processor). Note that the knowledge of $F$ is not
necessary to perform the comparisons in Step 3, since all the estimated
costs, especially $T_{\mathrm{MR}}(b, \lambda/2)$ for $\lambda/2$ Miller
Rabin tests, and 
\begin{multline*}
  T_{\mathrm{Galois}}(b, r, \dc, \dk) \simeq \\T_{\mathrm{MR}}(b, r) +
  T_{\mathrm{F}}(b, \dc) + T_{\mathrm{\zeta}}(b, \dc) +  T_{\sigma}(b, \dc,
  \dk) + T_{\mathrm{power}}(b, \dc, \dk)\,,
\end{multline*}
for Galois tests, are known up to $F$. Our
conclusions should thus be valid on any computer.

\begin{figure}[htb]
  \centering
  \includegraphics[scale=0.80]{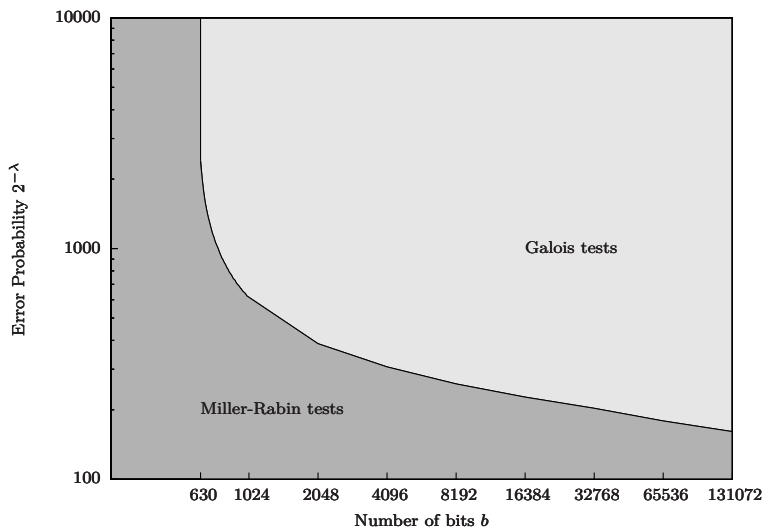}%
  \caption{Ranges of efficiency for the Galois test}
  \label{fig:galoistest}
\end{figure}

The  set of pairs  $(b,\lambda)$ for which  a Galois test is more
efficient than $\lambda/2$ Miller-Rabin tests is  the
pale domain 
in Figure~\ref{fig:galoistest}.
We observe
that when $b$ tends to  infinity, then the value of
$\lambda$ where the two 
methods cross tends to $47$.

A reasonably optimized implementation  in \textsc{magma v2.18-2} 
is available on the authors' web pages  for independent checks. 
In order to see how practical is this
implementation, we have picked a few
 random integers of sizes ranging from 1024 to
8192 bits, and we have measured the timings for those which turn to be
pseudo-primes. As expected, the cost ratio between
 $\lambda/2$ Miller-Rabin
tests and one equivalent  Galois test increases with $b$. Results are collected
in 
Table~\ref{tab:timings}.

\begin{table}[htb]
  \centering
  \begin{footnotesize}
  \begin{displaymath}
  \begin{array}{|c|c|c|c|c||r|r|r|r|r|r|r||r|}
    \hline
    \multicolumn{5}{|c||}{\mathrm{Parameters}} &
    \multicolumn{7}{|c||}{\mathrm{Galois}} &
    \mathrm{Miller-}\\
    \cline{1-12}
    b  & \lambda & r & \dc & \dk & 2 & 4 & 5 & 9 & \sigma(x) & x^n & \mathrm{Tot.} &
    \mathrm{Rabin} \\[0.1cm]\hline\hline
    \multirow{2}{*}{1024} & \multirow{2}{*}{512} & 129 & 1 & 15 & 0.0 & 0.3 & - & 0.0 & 0.0 & 0.2 & 0.5 & \multirow{2}{*}{0.5}\\
    & & 171 & 2 & 6 & 0.0 & 0.4 & 0.0 & 0.0 & 0.0 & 0.3 & 0.7 & \\\hline\hline
    \multirow{2}{*}{2048} & \multirow{2}{*}{1024} & 181 & 1 & 20 & 0.0 & 2.1 & - & 0.0 & 0.2 &
    1.5 & 3.8 & \multirow{2}{*}{6.2} \\
    & & 237 & 2 & 8 & 0.0 & 2.9 & 0.0 & 0.1 & 0.2 & 2.0 & 5.2 & \\\hline\hline
    \multirow{2}{*}{4096} & \multirow{2}{*}{2048} & 246 & 1 & 28 & 0.0 & 18.7 & - & 0.0 & 1.3 &
    11.6 & 31.6 & \multirow{2}{*}{75.1}\\
    & & 293 & 2 & 12 & 0.0 & 22.6 & 0.0 & 1.0 & 1.9 &
    19.8 & 45.3 & \\\hline\hline
    \multirow{3}{*}{8192} & \multirow{3}{*}{4096}
    & 333 & 1 & 40 & 0.4 & 176.9 & - & 0.7 & 1.3 & 122.5 & 314.8 & \multirow{3}{*}{1215.0}\\
     && 424 & 2 & 16 & 0.4 & 251.6 & 0.2 & 5.3 & 18.8 & 198.3 & 474.6 & \\
    & & 316 & 3 & 14 & 0.4 & 169.7 & 1.9 & 7.8 & 25.6 & 266.7  & 472.1 & \\\hline
  \end{array}
\end{displaymath}
\end{footnotesize}
\label{tab:timings}
\caption{Compared timings  for $b$-bit integers, and 
  prob. up to $2^{-{b}/{2}}$ (in seconds)}
\end{table}

\bibliographystyle{abbrv}

\end{document}

%% file: courbes.pspdftex
\begin{picture}(0,0)%
\includegraphics{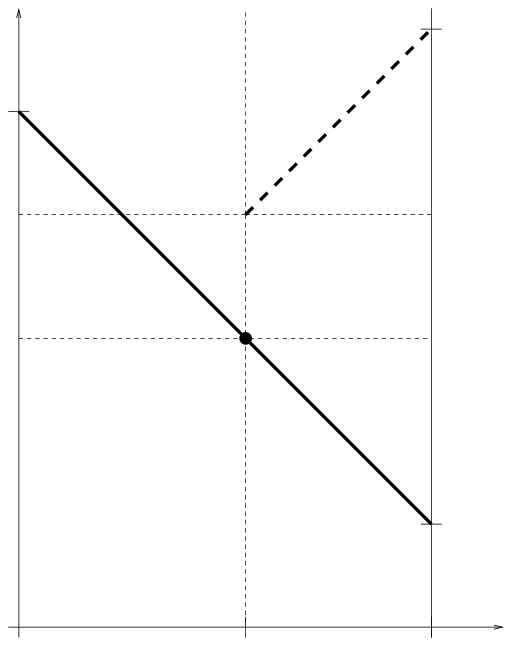}%
\end{picture}%
\setlength{\unitlength}{1302sp}%
\begingroup\makeatletter\ifx\SetFigFont\undefined%
\gdef\SetFigFont#1#2#3#4#5{%
  \reset@font\fontsize{#1}{#2pt}%
  \fontfamily{#3}\fontseries{#4}\fontshape{#5}%
  \selectfont}%
\fi\endgroup%
\begin{picture}(12705,9708)(661,-8257)
\put(4351,-1561){\makebox(0,0)[lb]{\smash{{\SetFigFont{7}{8.4}{\rmdefault}{\mddefault}{\updefault}{\color[rgb]{0,0,0}$1$}%
}}}}
\put(4201,-61){\makebox(0,0)[lb]{\smash{{\SetFigFont{7}{8.4}{\rmdefault}{\mddefault}{\updefault}{\color[rgb]{0,0,0}$1+\rho$}%
}}}}
\put(4051,-3361){\makebox(0,0)[lb]{\smash{{\SetFigFont{7}{8.4}{\rmdefault}{\mddefault}{\updefault}{\color[rgb]{0,0,0}$\delta+\rho$}%
}}}}
\put(5251,-8161){\makebox(0,0)[lb]{\smash{{\SetFigFont{7}{8.4}{\rmdefault}{\mddefault}{\updefault}{\color[rgb]{0,0,0}$0$}%
}}}}
\put(11251,-8161){\makebox(0,0)[lb]{\smash{{\SetFigFont{7}{8.4}{\rmdefault}{\mddefault}{\updefault}{\color[rgb]{0,0,0}$1$}%
}}}}
\put(8251,-8161){\makebox(0,0)[lb]{\smash{{\SetFigFont{7}{8.4}{\rmdefault}{\mddefault}{\updefault}{\color[rgb]{0,0,0}$1-\delta$}%
}}}}
\put(11851,-6061){\makebox(0,0)[lb]{\smash{{\SetFigFont{7}{8.4}{\rmdefault}{\mddefault}{\updefault}{\color[rgb]{0,0,0}$\rho$}%
}}}}
\put(11851,989){\makebox(0,0)[lb]{\smash{{\SetFigFont{7}{8.4}{\rmdefault}{\mddefault}{\updefault}{\color[rgb]{0,0,0}$1+\delta$}%
}}}}
\put(13351,-7636){\makebox(0,0)[lb]{\smash{{\SetFigFont{7}{8.4}{\rmdefault}{\mddefault}{\updefault}$\pi=\lim_{n\rightarrow\infty}\frac{\log\log\max\, p^v}{\log\log  n}$}}}}
\put(676,989){\makebox(0,0)[lb]{\smash{{\SetFigFont{7}{8.4}{\rmdefault}{\mddefault}{\updefault}{\color[rgb]{0,0,0}$\lim_{n\rightarrow \infty}-\frac{\log\log \mu}{\log\log  n}$}%
}}}}
\end{picture}%